\documentclass{article}
\usepackage{amsmath,amsfonts,amssymb,amscd}
\usepackage{mathptmx}     
\usepackage{latexsym}

\usepackage{graphicx}

\title{On Motzkin-Straus Type of  Results and Frankl-F\"uredi Conjecture for Hypergraphs}

\author{ Yuejian Peng \thanks{ College of Mathematics and Econometrics, Hunan University, Changsha 410082, P.R. China.  Email: ypeng1@163.com. Supported in part by National Natural Science Foundation of China (No. 11271116).}\and Yuping Yao \thanks{Corresponding author, College of Mathematics and Econometrics, Hunan University, Changsha 410082, P.R. China.Email: yupingyao1989@163.com}   }
\date{}

\date{}

\newtheorem{definition}{Definition}[section]
\newtheorem{theorem}{Theorem}[section]

\newtheorem{remark}[theorem]{Remark}

\newtheorem{lemma}[theorem]{Lemma}
\newtheorem{claim}[theorem]{Claim}

\newtheorem{con}[theorem]{Conjecture}

\newtheorem{prob}[theorem]{Problem}

\newcommand{\qed}{\hspace*{\fill} \rule{7pt}{7pt}}

\begin{document}
\maketitle

\begin{abstract}
A remarkable connection between the order of a maximum clique and
the Graph-Lagrangian of a graph was established by Motzkin and Straus in
1965. This connection and its extension were useful in both combinatorics and optimization. Since then, Graph-Lagrangian has been a useful tool in extremal combinatorics.  In this paper, we  give a  parametrized Graph-Lagrangian for non-uniform hypergraphs and provide  several  Motzkin-Straus type results for nonuniform hypergraphs which generalize results from \cite{PPTZ} and \cite{GLPS}. Another part of the paper concerns a long-standing conjecture of  Frankl-F\"uredi on Graph-Lagrangians of hypergraphs. We show the connection between the Graph-Lagrangian of $\{1, r_1, r_2, \cdots, r_l\}$-hypergraphs and $\{ r_1, r_2, \cdots, r_l\}$-hypergraphs.
Some of our results provide solutions to the maximum value of a class of polynomial functions over the standard simplex of the Euclidean space.

Keywords: Graph-Lagrangians of hypergraphs  Extremal combinatorics   Polynomial optimization
\textbf{\\AMS Subject Classification (2010)} 05C65, 05D05
\end{abstract}

\section{Introduction}
In 1965, Motzkin and Straus \cite{MS} established a connection between the order of a maximum clique and the Graph-Lagrangian of a graph.
This connection and its extensions were  successfully employed in optimization to provide heuristics for the
maximum clique problem \cite{B1,B2,B3,G9,PP}. This connection provided another proof of Tur\'an's
theorem \cite{Turan} which pushed the development of extremal graph theory.  More generally, the connection between Graph-Lagrangians and
Tur\'an densities can be used to give another proof of the
fundamental result of Erd\"os-Stone-Simonovits on Tur\'an densities
of graphs; see Keevash's survey paper \cite{Keevash}. However,  the obvious generalization of Motzkin and Straus' result to  $r$-uniform hypergraphs is false.  i.e., the Graph-Lagrangian of a hypergraph is not always the same as the Graph-Lagrangian of its maximum cliques. There are many examples of  $r$-uniform hypergraphs other than complete $r$-uniform hypergraphs  that do not achieve their Graph-Lagrangian on any proper subhypergraph. In spite of this,  Graph-Lagrangians has been a useful tool in extremal problems in combinatorics.
In 1980's, Sidorenko \cite{Sidorenko} and Frankl and F\"uredi \cite{FF}
developed the method of applying Graph-Lagrangians in determining
hypergraph Tur\'an densities. More recent  applications of Graph-Lagrangians can
be found in Keevash's survey paper (\cite{Keevash}),  \cite{Mu} and \cite{keevash2}.   In most applications in extremal combinatorics, we need an upper bound for the Graph-Lagrangians of hypergraphs.
In the course of estimating Tur\'an densities of hypergraphs by applying the Graph-Lagrangians of related hypergraphs, Frankl and F\"uredi \cite{FF} asked the following question: Given $r \ge 3$ and $m \in {\mathbb N}$ how large can the Graph-Lagrangian of an $r$-graph with $m$ edges be? They proposed the following conjecture: The $r$-graph with $m$ edges formed by taking the first $m$ sets in the colex ordering of $N^{(r)}$ has the largest Graph-Lagrangian of all $r$-graphs with $m$ edges.  Moztkin-Straus result implies that this conjecture is true for $r=2$. For $r\ge 3$, this conjecture seems to be very challenging. Talbot first confirmed this conjecture for some cases  in \cite{talbot}. Later  Tang et al. confirmed this conjecture for some more cases in  \cite{TPZZ,TPZZ1,STZP}.


  Recently, the study of Tur\'an
densities of non-uniform hypergraphs has been motivated by the study
of extremal poset problems \cite{GK,GL}. In \cite{JL}, Johnston and Lu gave a generalization of the concept of Tur\'an density of a
non-uniform hypergraph. In \cite{PPTZ},  Peng et al. introduced the Graph-Lagrangian of a non-uniform
hypergraph, and gave an extension of
Erd\"os-Stone-Simonovits theorem to non-uniform hypergraphs whose edges contain $1$  or
$2$ vertices by applying Graph-Lagrangians of non-uniform hypergraphs (this extension of Erd\"os-Stone-Simonovits theorem to non-uniform hypergraphs was given in \cite{JL} by a different method). In this paper, we study a more generalized  question for non-uniform hypergraphs and provide several results related to this question( Theorems \ref{th2}, \ref{thm4},  \ref{th4} and \ref{1.2}).  Although the truth of Conjecture \ref{conjecture} of Frankl and F\"uredi  is not known in general even for $r$-uniform hypergraphs, we propose that a similar result is true for non-uniform hypergraphs (Problem 2) and provide some partial results (Theorem \ref{1connection}).

   Our main results provide solutions to the maximum value of a class of polynomial functions in several variables.


\section{Definitions, notations and main results}

 A \emph{hypergraph} is a pair $H = (V(H),E(H))$ consisting of a vertex set
$V(H)$ and an edge set  $E(H)$, where each edge is a
subset of $V(H)$. The set
$T(H)=\{|e|:e \in E\}$ is called the set of \emph{edge types} of
$H$. We also say that $H$ is a $T(H)$-graph. For example, if
$T(H)=\{1,3\}$, then we say that $H$ is a $\{1,3\}$-graph. If
all edges have the same cardinality $r$, then $H$ is an $r$-uniform
hypergraph, which is simply written as $r$-graph. A $2$-uniform
hypergraph is  a simple graph. A hypergraph is non-uniform
if it has at least two edge types. Write $H_n^T$ for
a hypergraph $H$ on $n$ vertices with $T(H)=T$. For any $r \in T(H)$, the
\emph{rth-level hypergraph} $H^r$ is the hypergraph consisting of all
edges containing $r$ vertices of $H$. For $Q\subset T$, let $H^Q$ denote the hypergraph $\cup_{r\in Q} H^r$.
We also use  $E^r$ to denote
the set of all edges with $r$ vertices of $H$.  For convenience, an
edge $\{i_1, i_2, \ldots, i_r\}$ in a hypergraph is simply written
as $i_1 i_2 \ldots i_r$ throughout the paper.

For an integer $n$, let $[n]$ denote the set $\{1, 2, \cdots, n\}$. For a set $V$ and a positive integer $i$, let ${V \choose i}$ be the set of all subsets of $V$ with $i$ elements. The complete hypergraph $K_n^T$ is a hypergraph on  vertex set $[n]$  with edge set
$\bigcup\limits_{i \in T} {\left( {\begin{array}{*{20}{c}}
{[n]}\\
i
\end{array}} \right)}$. For example, $K_n^{\{ r\} }$ is the complete $r$-uniform
hypergraph on $n$ vertices.
$K_n^{[r]}$ is the non-uniform hypergraph with all possible edges of cardinality at
most $r$. Let $[n]^{T}$ represent the complete $T$-type hypergraph
on vertex set $[n]$. For example, $[n]^{\{1, 3\}}$ represents the complete $\{1, 3\}$-hypergraph
on vertex set $[n]$. We also let $[n]^{(r)}$ represent the complete $r$-uniform hypergraph
on vertex set $[n]$.
 A hypergraph  $H$ is  a {\it subgraph} of a hypergraph $G$, denoted by $H\subseteq G$ if $V(H)\subseteq V(G)$ and $E(H)\subseteq E(G)$. A complete subhypergraph of a hypergraph $H$ with the same edge type as $T(H)$ is called a clique of $H$.  If $W\subseteq V(H)$, then the subhypergraph of $H$ induced by $W$ is denoted by $H[W]$, i.e. the vertex set of $H[W]$ is $W$ and the edge set of $H[W]$ is the set of all edges in $H$ whose vertices are in $W$.

\begin{definition} \label{definitionGraph-Lagrangian}
For  an $r$-uniform graph $H$ with the vertex set $[n]$, edge set $E(H)$, and a vector $\vec{x}=(x_1,\ldots,x_n) \in {\mathbb R}^n$, associate a homogeneous polynomial in $n$ variables, denoted by  $\lambda (H,\vec{x})$  as follows:
$$\lambda (H,\vec{x}):=\sum_{i_1i_2 \cdots i_r \in E(H)}x_{i_1}x_{i_2}\ldots x_{i_r}.$$
Let $S:=\{\vec{x}=(x_1,x_2,\ldots ,x_n): \sum_{i=1}^{n} x_i =1, x_i
\ge 0 {\rm \ for \ } i=1,2,\ldots , n \}$.
The Graph-Lagrangian of $H$, denoted by $\lambda (H)$, is the maximum
 of the above homogeneous  multilinear polynomial of degree $r$ over the standard simplex $S$. Precisely,
$$\lambda (H): = \max \{\lambda (H, \vec{x}): \vec{x} \in S \}.$$
\end{definition}
The value $x_i$ is called the {\em weight} of the vertex $i$.
A vector $\vec{x}=(x_1, x_2, \ldots, x_n) \in {\mathbb R}^n$ is called a feasible weighting for $H$ if and only if
$\vec{x}\in S$. A vector $\vec{y}\in S$ is called an {\em optimal weighting} for $H$
if and only if $\lambda (H, \vec{y})=\lambda(H)$.

\begin{remark}
$\lambda (H)$ was called Lagrangian of $H$ in literature \cite{FF,talbot,FR,PZ}. The terminology `Graph-Lagrangian' was suggested by Franco Giannessi.
\end{remark}

Motzkin and Straus in \cite{MS} proved the following result for the
Graph-Lagrangian of a $2$-graph. It shows that the Graph-Lagrangian of a graph
is determined by the order of its maximum cliques.

\begin{theorem}\label{MStheo}\cite{MS}
If G is a $2$-graph in which a largest clique has order $t$, then,
\[\lambda (G) = \lambda \left( {{K_t}^{\{ 2\} }} \right) = \lambda
\left( {{{\left[ t \right]}^{(2)}}} \right) = \frac{1}{2}\left( {1 - \frac{1}{t}} \right).\]
\end{theorem}

\noindent  The Motzkin-Straus result and its extension had many applications in extremal problems in graphs and hypergraphs \cite{Keevash}.   However,  the obvious generalization of Motzkin and Straus' result to  $r$-uniform hypergraphs is false.  i.e., the Graph-Lagrangian of a hypergraph is not always the same as the Graph-Lagrangian of its maximum cliques. 
In spite of this, there are still applications of Graph-Lagrangians of hypergraphs in  determining hypergraph Tur\'an densities \cite{Keevash,Sidorenko,FF1}.
 In most applications, we need an upper bound for the Graph-Lagrangians of hypergraphs.
 Frankl and F\"uredi \cite{FF} asked the following question: Given $r \ge 3$ and $m \in {\mathbb N}$ how large can the Lagrangian of an $r$-graph with $m$ edges be? In order to state their conjecture on this problem we require the following definition. For distinct $A, B \in {\mathbb N}^{(r)}$ we say that $A$ is less than $B$ in the {\em colex ordering} if $max(A \triangle B) \in B$, where $A \triangle B=(A \setminus B)\cup (B \setminus A)$. For example we have $246 < 156$ in ${\mathbb N}^{(3)}$ since $max(\{2,4,6\} \triangle \{1,5,6\}) =\{5\}\in \{1,5,6\}$. In Colex ordering, $123<124<134<234<125<135<235<145<245<345<126<136<236<146<246<346<156<256<356<456<127<\cdots .$
Let $C_{m,r}$ denote the $r$-graph with $m$ edges formed by taking the first $m$ elements in the colex ordering of ${\mathbb N}^{(r)}$. When $m={t \choose r}$,  the $r$-graph  $C_{ {t \choose r},r}$  is $[t]^{(r)}$. The following conjecture of Frankl and F\"uredi (if it is true) proposes a  solution to the above question.

\begin{con} (Frankl and F\"uredi \cite{FF})\label{conjecture} The $r$-graph formed by taking the first $m$ sets in the colex ordering of $N^{(r)}$ has the largest Lagrangian of all $r$-graphs with $m$ edges.  In other words, if $H$ is an $r$-graph with $m$ edges, then $\lambda(H)\le \lambda(C_{m, r})$.
\end{con}

Motzkin-Straus's Theorem (Theorem \ref{MStheo}) implies that this conjecture is true when $r=2$ by Theorem \ref{MStheo}. For the case $r=3$,  Talbot in \cite{talbot} proved the following.

\begin{theorem} (Talbot \cite{talbot}) \label{Tal} Let $m$ and $t$ be integers satisfying
$${t \choose 3}-2 \le m \le {t \choose 3} + {t-1 \choose 2} -t.$$
Then Conjecture \ref{conjecture} is true for $r=3$ and this value of $m$.
\end{theorem}

Recently, Tang et al. verified this conjecture for more cases.

\begin{theorem}  \cite{talbot,TPZZ,STZP} \label{TPZZ} Let $m$ and $t$ be integers satisfying
$${t \choose 3}-7 \le m \le {t \choose 3} + {t-1 \choose 2} -\frac{1}{2}t.$$
Then Conjecture \ref{conjecture} is true for $r=3$ and this value of $m$.
\end{theorem}

Let
$$\lambda_{(m, n)}^r=max \{\lambda(H): \,\ H=(V,E)\,\ {\rm is \,\ an \,\ r-graph}, |V|=n,|E|=m\}.$$
For $r\ge 4$, the only known results are



\begin{theorem} (Talbot \cite{talbot}) \label{Talr}For any $r \ge 4$ there exists constants $\gamma_r$ and $\kappa_0(r)$ such that if $m$ satisfies
$${t \choose r} \le m \le {t \choose r} + {t-1 \choose r-1} - \gamma_r (t)^{r-2},$$
with $t \ge \kappa_0(r)$, then $\lambda_{(m, t+1)}^r=\lambda(C_{m,r})=\lambda([t]^{(r)})$.
\end{theorem}

\begin{theorem}  \cite{TPZZ1} \label{TPZZ1} Let $m$, $r$ and $t$ be integers satisfying
$${t \choose r}-4 \le m \le {t \choose r}.$$
Then $\lambda_{(m, t)}^r=\lambda(C_{m,r})$.
\end{theorem}

Recently, the study of Tur\'an
densities of non-uniform hypergraphs has been motivated by the study
of extremal poset problems \cite{GK,GL}. In
\cite{JL},  Johnston and Lu gave a generalization of the
concept of Tur\'an density to a non-uniform hypergraph. In \cite{PPTZ}, Peng et al. generalized the concept of Graph-Lagrangian to non-uniform hyergraphs, gave a generalization of Mozkin-Straus result to $\{1,2\}$-graphs, and consequently applied it obtaining a result on Tur\'an densities of $\{1, 2\}$-graphs  similar to Erd\H os-Stone-Simonovits classical result on Tur\'an densities of graphs.
In this paper, we study the following general optimization problem for
 non-uniform hypergraphs which generalizes the concept of Graph-Lagrangians.

\begin{prob}\label{generalproblem}
Let $H$ be  an $\{r_0,r_1,r_2,\ldots,r_l\}$-graph, $r_0<r_1<r_2< \ldots <r_l$, with vertex set $V(H)=[n]$ and
edge set $E(H)$.  Let $S=\{\vec{x}=(x_1,x_2,\ldots ,x_n)\in \mathbb{R}^n: \sum_{i=1}^{n} x_i =1, x_i
\ge 0 {\rm \ for \ } i=1,2,\ldots , n \}$. Let $\alpha_1,\alpha_2,\cdots,\alpha_l$ be non-negative constants. For $\vec{x}\in S$,
let
\begin{eqnarray*}
L_{\{\alpha_1,\alpha_2,\cdots,\alpha_l\}} (H,\vec{x})&:=&\sum_{ i_1i_2 \ldots i_{r_0}  \in E(H^{r_0})}x_{i_1}x_{i_2}\ldots x_{i_{r_0}}+\alpha_1 \sum_{ i_1i_2 \ldots i_{r_1}  \in E(H^{r_1})}x_{i_1}x_{i_2}\ldots x_{i_{r_1}}\\
&+&\ldots+\alpha_{l} \sum_{ i_1 i_2\ldots i_{r_l} \in E(H^{r_l})}x_{i_1}x_{i_2}\ldots x_{i_{r_l}}.
\end{eqnarray*}
The {\em polynomial optimization problem} of $H$ is
\begin{eqnarray}\label{lg}
L_{\{\alpha_1,\alpha_2,\cdots,\alpha_l\}}  (H): = \max \{L(H, \vec{x}): \vec{x} \in S \}.
\end{eqnarray}
We sometimes simply write $L_{\{\alpha_1,\alpha_2,\cdots,\alpha_l\}} (H,\vec{x})$ and $L_{\{\alpha_1,\alpha_2,\cdots,\alpha_l\}}(H)$ as $L(H,\vec x)$ and $L(H)$ if there is no confusion.
The value $x_i$ is called the {\em weight} of the vertex $i$. A vector $\vec{x}=(x_1, x_2, \ldots, x_n) \in {\mathbb R}^n$ is called a feasible solution to (\ref{lg}) if and only if $\vec{x}\in S$. A vector $\vec{y}\in S$ is called a {\em solution} to optimization problem (\ref{lg}) if and only if $L (H, \vec{y})=L(H)$.
\end{prob}

\begin{remark} \label{re_1}
If $G$ is a subhypergraph of $H$, then $L_{\{\alpha_1,\alpha_2,\cdots,\alpha_l\}} (G) \le L_{\{\alpha_1,\alpha_2,\cdots,\alpha_l\}} (H)$.
\end{remark}

The \emph{characteristic vector} of a set $U$,   denoted by $\vec{x}^{U}=(x_1^{U}, x_2^{U}, \ldots, x_n^{U})$,  is the vector in $S$ defined as: $$x_i^{U}=\frac{1_{i\in U}}{|U|}$$ where $|U|$ denotes the cardinality of $U$ and $1_{P}$ is the indicator function returning 1 if property $P$ is satisfied and 0 otherwise.

In this paper, we show the following  result to Problem \ref{generalproblem} for $\{1, 2\}$-graphs which generalizes a result in \cite{PPTZ}.

\begin{theorem} \label{th2}
Let $\alpha_2>0$ be a constant. If $H$ is a $\{1,2\}$-graph with $n$ vertices and the order of its maximum clique is $t$, where $t\ge \alpha_2$, then $ L_{\{ \alpha_2 \}}(H)= L_{\{ \alpha_2 \}}(K^{\{1,2\}}_t)=1+\frac{\alpha_2}{2} - \frac{\alpha_2}{2t}$. Furthermore, the characteristic vector of a maximum clique is a solution to optimization problem (\ref{lg}).
\end{theorem}

In \cite{GLPS}, Gu et al. give some Motzkin-Straus type results to non-uniform hypergraphs. In a similar way, we give Motzkin-Straus type results to $\{1,r \}$-graphs and  $\{1, 2, 3\}$-graphs regarding Problem (\ref{lg}).

\begin{theorem}\label{thm4}
Let $\alpha_1>0$ be a constant. Let $H$ be a $\{1,r\}$-graph. If both the order of its maximum complete $\{1,r\}$-subgraphs and  the order of its maximum complete $\{1\}$-subgraphs are $t$, where $\displaystyle {t\geq \lceil {[\alpha_1-(r-2)!]^{r-2} \over (r-2)!\alpha_1^{r-3}}\rceil }$, then \[L _{\{ \alpha_1 \}}(H) = L_{\{ \alpha_1 \}} \left( {{K_t}^{\{ 1,r\} }} \right) ={1+ \alpha_1 \frac{\prod_{i=1}^{r-1} (t-i)}{r!t^{r-1}}}.\]
Furthermore, the characteristic vector of a maximum clique is a solution to optimization problem (\ref{lg}).
\end{theorem}

\begin{theorem}\label{th4}
Let $\alpha_1,\alpha_2>0$ be constants. Let $H$ be a $\{1, 2, 3\}$-graph. If both the order of its maximum complete
$\{1,2, 3\}$-subgraph  and  the order of its maximum complete
$\{1\}$-subgraph are $t$, where $t\ge \lceil{ (\alpha_1+\alpha_2)^2-\alpha_2 \over \alpha_1+\alpha_2}\rceil $, then
\[ L_{\{ \alpha_1,\alpha_2 \}}(H) =  L_{\{ \alpha_1, \alpha_2 \}} \left( {{K_t}^{\{ 1, 2, 3\} }} \right) ={1+  \alpha_1 \frac{t-1}{2t}+  \alpha_2 \frac{(t-1)(t-2)}{6t^2}}.\]
Furthermore, the characteristic vector of a maximum clique is a solution to optimization problem (\ref{lg}).
\end{theorem}

A result to $\{1,r_1,r_2,\ldots,r_l\}$- graph will be also given.
\begin{definition} \label{defi1.1}
Let $H$ be  an $\{1,r_1,r_2,\ldots,r_l\}$-graph, $1<r_1<r_2< \ldots <r_l$, with vertex set $V(H)=[n]$ and
edge set $E(H)$. For $i \in V(H^1) $, $i$ is isolated in $H$ if and only if there is no edge $e \in E(H^{r_i})$ $(i=1, \ldots ,l)$ such that $i\in e$. The set of all isolated vertices of $H$ is denoted by $D(H)$.
\end{definition}
For example, if $H=\{1,2,3,4,5\} \cup \{12,13\} \cup \{123,356\} $, then vertex $4$ is isolated in $H$ and vertices $4$ and $5$ are isolated in $H[V(H^1)]$, so $D(H[V(H^1)])=\{4,5\}$.

\begin{theorem}\label{1.2}
Let $\alpha_1,\alpha_2,\ldots ,\alpha_l$ be non-negative constants. Let H be a $\{1,r_1,r_2,\ldots,r_l\}$- graph, where $1<r_1<\ldots<r_l$ and $\sum_{i=1,\ldots,l}{\alpha_i \over (r_i-1)!} \le 1$. If $V(H^1) \not =  D(H[V(H^1)])$, then $ L_{\{\alpha_1,\alpha_2,\cdots,\alpha_l \}}(H)= L_{\{\alpha_1,\alpha_2,\cdots,\alpha_l \}}(H[V(H^1)\setminus D(H[V(H^1)])])$; and if $V(H^1)=D(H[V(H^1)])$, then $ L_{\{\alpha_1,\alpha_2,\cdots,\alpha_l \}}(H)=1$.
\end{theorem}

Although the truth of Conjecture \ref{conjecture} is not known in general even for $3$-uniform hypergraphs, we propose a similar question for hypergraphs.  Similarly, for distinct  sets $A, B \subset {\mathbb N}$,  we say that $A$ is less than $B$ in the {\em colex ordering} if $max(A \triangle B) \in B$.  Let $C_{m, T}$ denote the hypergraph with edge type $T$ and $m$ edges formed by taking the first $m$ elements in the colex ordering.

\begin{prob}\label{generalff} Let $H$ be a hypergraph with edge type $T=\{r_0,r_1,r_2,\ldots,r_l\}$ and $m$ edges.  For  what conditions on $\alpha_i>0$, the inequality
\begin{eqnarray*}
L_{\{\alpha_1,\alpha_2,\cdots,\alpha_l\}}  (H) \le L_{\{\alpha_1,\alpha_2,\cdots,\alpha_l\}}  (C_{m, T})
\end{eqnarray*}
holds?
\end{prob}

Theorems \ref{MStheo} and \ref{th2} provided some results to Problem \ref{generalff}  for $T=\{2\}$ or  $T=\{1, 2\}$. For $T=\{3\}$, Theorems \ref{Tal} and \ref{TPZZ} provided partial results to this problem. 
We show the following connection between $\{1,r_1,r_2,\ldots,r_l\}$-hypergraphs and $\{r_1,r_2,\ldots,r_l\}$-hypergraph concerning this question.


\begin{theorem}\label{1connection}
Let $r_1,\ldots,r_l$ be positive integers satisfying $1<r_1<\ldots<r_l$. Let $\alpha_i(i=1,\ldots,l)$ be positive constants satisfying $\sum_{i=1}^l {\alpha_i \over (r_i-1)!} \le 1$. Let $m$ and $t$ be positive integers satisfying $t+\sum_{i=1}^l {t \choose r_i}<m\le t+1+\sum_{i=1}^l {t+1 \choose r_i}$.  Let H be a $\{1,r_1,r_2,\ldots,r_l\}$-hypergraph with $m$ edges and $n$ vertices. If for an $\{r_1,r_2,\ldots,r_l\}$-hypergraph $G$ with $m-t-1$ edges and $n$ vertices,
$L_{\{\alpha_2,\cdots,\alpha_l\}}  (G) \le L_{\{\alpha_2,\cdots,\alpha_l\}}  (C_{m-t-1, \{r_1,r_2,\ldots,r_l\}})$ holds, then
  $L_{\{\alpha_1,\alpha_2,\cdots,\alpha_l\}}  (H) \le L_{\{\alpha_1,\alpha_2,\cdots,\alpha_l\}}  (C_{m, \{1,r_1,r_2,\ldots,r_l\}})$ holds.
\end{theorem}

Combining this theorem and the known results given in Theorem \ref{TPZZ}, we can get corresponding results for $\{1, 3\}$-hypergraphs.

\section{Some preliminaries}

In this section, we give some preliminary results to be applied in the proof.

The {\em support} of a vector $\vec{x}\in S$, denoted by $\sigma(\vec{x})$, is the set of indices corresponding to positive components of $\vec{x}$, i.e.,
$$\sigma(\vec{x})=\{i: x_i>0, 1\le i\le n\}.$$
We will impose an additional condition on a solution  $\vec x = (x_1,\cdots, x_n)$  to optimization problem (\ref{lg}).

(*) $\vert\sigma(x)\vert$ is minimal, i.e., if $\vec y$ is a legal weighting for $H$ satisfying $\vert\sigma(y)\vert<\vert\sigma(x)\vert$, then $L \left( {{H},\vec y} \right)<L \left( H \right)$.

For a hypergraph $H=(V,E)$, $i\in V$, and $r\in T(H)$, let $E_i^r=\{A \in V^{(r-1)}: A \cup \{i\} \in E^r\}$. For a pair of vertices $i,j \in V$, let $E_{ij}^r=\{B \in V^{(r-2)}: B \cup \{i,j\} \in E^r\}$.
Let
$(E^{r}_i)^c=\{A \in V^{(r-1)}: A \cup \{i\} \in V^{(r)} \backslash E\}$, $(E^r_{ij})^c=\{B \in V^{(r-2)}: B \cup \{i,j\} \in V^{(r)} \backslash E^r\}$, and $E_{i\setminus j}^r=E_i^r\cap (E^{r}_{j})^c.$
Let $L(E_i^r, \vec{x})=\alpha_{r}\lambda (E_i^r, \vec{x})$, $L(E_{ij}^r, \vec{x})=\alpha_{r}\lambda (E_{ij}^r, \vec{x})$, and $L(E_{i\setminus j}^r, \vec{x})=\alpha_{r}\lambda (E_{i\setminus j}^r, \vec{x})$, where $\alpha_{r_0}=1$. Let $E_i=\cup_{r\in T(H)} E_i^r$, $E_{i\setminus j}=\cup_{r\in T(H)} E_{i\setminus j}^r$, and   $E_{ij}=\cup_{r\in T(H)} E_{ij}^r$. Let  $L(E_i, \vec{x})=\sum_{r\in T(H)} L(E_i^r, \vec{x})$,  $L(E_{ij}, \vec{x})=\sum_{r\in T(H)} L(E_{ij}^r, \vec{x})$, and  $L(E_{i\setminus j}, \vec{x})=\sum_{r\in T(H)} L(E_{i\setminus j}^r, \vec{x})$. Note that $L(E_i, \vec{x})={\partial L (H, {\vec x})\over \partial x_i}$ and $L(E_{ij}, \vec{x})={\partial L (H, {\vec x})\over \partial x_i\partial x_j}$.


\begin{lemma}\label{FR}\cite{FR}
Let $H = (V, E)$ be an $r$-graph and $\vec x = (x_1,\cdots, x_n)$ be
an optimal legal weighting for $H$ satisfying (*).  Then for  $i, j \in \sigma(\vec{x})$,

(a) $\lambda(E_i , \vec x) =\lambda(E _j, \vec x) =
r\lambda(G)$,

(b) there is an edge in $E$ containing both $i$ and
$j$.
\end{lemma}

 We give a similar result
for a non-uniform hypergraph below.

\begin{lemma} \label{Lemmkkt}
Let ${\vec x}=(x_1, x_2, \ldots, x_n)$ be a solution to  the polynomial programming (\ref{generalproblem}) satisfying (*). Then for  $i, j \in \sigma(\vec{x})$,

(a) ${\partial L (H, {\vec x})\over \partial x_i}={\partial  L (H, {\vec x})\over \partial x_j}.$  This is equivalent to  $L(E_i, {\vec x})=L(E_j, {\vec x})$,

(b) there exists an edge $e \in E(H)$  such that $\{i,j\} \subseteq e$.
\end{lemma}

{\em Proof.} (a) Suppose, for a contradiction, that there exist $i, j\in \sigma(\vec{x})$ such that ${\partial L (H, {\vec x})\over \partial x_i}> {\partial  L (H, {\vec x})\over \partial x_j}$. We define a new feasible solution $\vec y$ to $(\ref{lg})$ as follows. Let $y_q=x_q$ for $q \neq i,j$, $y_i=x_i+ \delta$ and $y_j=x_j- \delta \geq 0,$ then
\begin{eqnarray*}
&& L(H, \vec y)-  L(H, \vec x)\\
&=&\delta({\partial  L (H, {\vec x})\over \partial x_i}- x_j{\partial^2  L (H, {\vec x})\over {\partial x_i \partial x_j}})-
\delta({\partial  L (H, {\vec x})\over \partial x_j}- x_i{\partial^2  L (H, {\vec x})\over {\partial x_i \partial x_j}})\\
& &+(\delta x_j-\delta x_i-\delta^2){\partial^2  L (H, {\vec x})\over {\partial x_i \partial x_j}}\\
&=&\delta({\partial  L (H, {\vec x})\over \partial x_i}- {\partial  L (H, {\vec x})\over \partial x_j})-\delta^2{\partial^2  L (H, {\vec x})\over {\partial x_i \partial x_j}}\\
&>&0
\end{eqnarray*}
for some small enough $\delta$, contradicting to that $\vec x$ is a solution to optimization problem (\ref{lg}).

(b) Suppose, for a contradiction, that there exist $i,j\in \sigma(\vec{x})  $ such that $\{i,j\} \nsubseteq e$ for any $e \in E(H)$. We define a new feasible solution $\vec y$ to $(\ref{lg})$ as follows. Let $y_q=x_q$ for $q \neq i,j$, $y_i=x_i+ x_j$ and $y_j=x_j- x_j=0,$ then $\vec y$ is clearly a feasible solution for $H$, and $$ L(H, \vec y)-  L(H, \vec x)=x_j({\partial  L (H, {\vec x})\over \partial x_i}- {\partial  L (H, {\vec x})\over \partial x_j})-x_j^2{\partial^2  L (H, {\vec x})\over {\partial x_i \partial x_j}}=0.$$
So $\vec y$ is a solution to optimization problem (\ref{lg}) and $\vert\sigma(y)\vert=\vert\sigma(x)\vert -1$, contradicting the minimality of $\vert\sigma(x)\vert$. \qed


In \cite{talbot},  Talbot introduced the definition of a left-compressed $r$-uniform hypergraph.  This concept is generalized to   non-uniform hypergraphs in \cite{GLPS}.

Let $H=([n],E)$ be a $T(H)$-graph, where $n$ is a positive integer. For $e \in E$, and $i,j\in [n]$ with
$i<j$, define

\begin{equation}
    {C_{i\leftarrow j}}\left( e \right)=
   \begin{cases}
   {(e\backslash \{ j\} ) \cup \{ i\} } &\mbox{if $i \notin e$ and $j\in e$,}\\
   e &\mbox{otherwise.}
   \end{cases}
  \end{equation}
And
\begin{equation}\label{Lij}
\mathcal{C}_{i\leftarrow j}(E)=\{C_{i\leftarrow j}(e):e \in E\}\cup \{e:e, C_{i\leftarrow j}\left( e \right) \in E\}.
\end{equation}
Note that $|\mathcal{C}_{i\leftarrow j}(E)|=|E|$ from the definition of $\mathcal{C}_{i\leftarrow j}(E)$.

We say that $E$ or $H$ is \emph{left-compressed} if and only if
$\mathcal{C}_{i\leftarrow j}(E)=E$ for every $1\leq i<j$. If a $T(H)$-hypergraph $H$ is left-compressed, then for every $r\in T(H)$, the
$r$-level hypergraph $H^r$ is left-compressed. An equivalent perhaps more intuitive
definition of left-compressed hypergraph is that a $T(H)$-hypergraph $H=([n],E)$ is
left-compressed if and only if for any $r \in T(H)$, $j_1j_2 \cdots j_r \in E$
implies $i_1i_2\cdots i_r \in E$ provided $i_p \le j_p$ for every $p$, $1 \le
p\le r$.
 Moreover, if
$H$ is a left-compressed $T(H)$-hypergraph and $i<j$, then for every $r\in T(H)$,
$E_{j\setminus i}^r=\emptyset$.

The following lemma is similar to a result given in \cite{talbot}.

\begin{lemma}\label{lem1}
Let $H=([n],E)$ be a $T(H)$-graph,
$i,j\in [n]$ with $i<j$ and $\vec x = (x_1,\cdots, x_n)$ be an optimal legal weighting of $H$. Write $H_{i\leftarrow j}=([n],\mathcal{C}_{i\leftarrow j}(E))$. Then,
$$L(H, \vec x)\leq L(H_{i\leftarrow j}, \vec x).$$
\end{lemma}

\begin{proof}
If $1 \notin T(H)$, then,
$$L(H_{i\leftarrow j}, \vec x)-L(H, \vec x)=\sum\limits_{r \in T(H)} {\sum\limits_{\scriptstyle e \in E^r,{C_{i\leftarrow j}}\left( e \right) \notin E^r\hfill\atop
\scriptstyle i \notin e,j \in e\hfill} {L(e\backslash \{ j\} ,\vec x)\left( {{x_i} - {x_j}} \right)} } ,$$
and if $1 \in T(H)$, then,
$$L(H_{i\leftarrow j}, \vec x)-L(H, \vec x)=
\sum\limits_{\scriptstyle r \in T(H)\hfill\atop
\scriptstyle r \ge 2\hfill} {\sum\limits_{\scriptstyle e \in {E^r},{C_{i\leftarrow j}}\left( e \right) \notin {E^r}\hfill\atop
\scriptstyle i \notin e,j \in e\hfill} {L(e\backslash \{ j\} ,\vec x)\left( {{x_i} - {x_j}} \right)} }  + \left( {{x_i} - {x_j}} \right)I,$$
where $I$ satisfies that $I=1$, if $i\notin E^1$
$j\in E^1$, and otherwise $I=0$.
Hence $L(H_{i\leftarrow j}, \vec x)-L(H, \vec x)$ is
nonnegative in any case, since $i<j$ implies that
$x_i\geq x_j$. So this lemma holds.
\end{proof}

\begin{remark}\label{r1} Let ${\vec x}=(x_1, x_2, \ldots, x_n)$ be  a solution to the optimization problem (\ref{lg}). Let $i, j  \in \sigma(x)$ with $i<j$.

 (a) Lemma \ref{Lemmkkt} part (a) implies that
$$x_jL(E_{ij}, {\vec x})+L (E_{i\setminus j}, {\vec x})=x_iL(E_{ij}, {\vec x})+L(E_{j\setminus i}, {\vec x}).$$
In particular, if $H$ is left-compressed, then
$$(x_i-x_j)L(E_{ij}, {\vec x})=L(E_{i\setminus j}, {\vec x})$$
since $E_{j\setminus i}=\emptyset$.

(b) If  $H$ is left-compressed, then
\begin{equation}\label{enbhd}
x_i-x_j={L (E_{i\setminus j}, {\vec x}) \over L(E_{ij}, {\vec x})}
\end{equation}
holds.  If  $H$ is left-compressed and  $E_{i\setminus j}=\emptyset$, then $x_i=x_j$.

(c)If  $H$ is left-compressed, then
\begin{equation}\label{conditiona}
x_1 \ge x_2 \ge \ldots \ge x_n \ge 0.
\end{equation}
\end{remark}

A result similar to Lemma $2.4$ in \cite{talbot} is also true for non-uniform hypergraphs.

\begin{lemma} \label{LemmaTal7nonu}
For any positive integers $m, t, r_1, r_2, \ldots, r_l$ satisfying 
$$\sum_{i=1}^l {t \choose r_i} \le m \le \sum_{i=1}^l {t \choose r_i}
\\+ \sum_{i=1}^l {t-1 \choose r_i-1},$$
we have $L(C_{m, \{r_1, r_2, \ldots, r_l\}}) = L([t]^{\{r_1, r_2, \ldots, r_l\}})$.
\end{lemma}
{\bf Proof.} Note that the vertex set of $C_{m, \{r_1, r_2, \ldots, r_l\}}$ is $[t+1]$ and $[t]^{\{r_1, r_2, \ldots, r_l\}}\subset C_{m,\{r_1, r_2, \ldots, r_l\}}$.  So $L(C_{m, \{r_1, r_2, \ldots, r_l\}}) \ge L([t]^{\{r_1, r_2, \ldots, r_l\}})$. Let ${\vec x}=(x_1, x_2, \ldots, x_{t+1})$ be an optimal weighting of $C_{m, \{r_1, r_2, \ldots, r_l\}}$. Note that the range of $m$ guarantees that there is no edge in $C_{m, \{r_1, r_2, \ldots, r_l\}}$ containing both $t$ and $t+1$. By Lemma \ref{Lemmkkt}, $x_{t+1}=0$.   Therefore, 
$$L(C_{m,\{r_1, r_2, \ldots, r_l\}})\le L([t]^{\{r_1, r_2, \ldots, r_l\}}).$$ 
This completes the proof of this lemma.  \qed

Denote $L_{(m, n)}^T=\max\{L(H): H $ is a  $T$-hypergraph with  $m$ edges and no  more  than $n$ vertices  \}.

\begin{definition}
We say that $H$ is an extremal $T$-graph for $L_{(m, n)}^T$ if $H$ is a $T$-graph  with $m$ edges and no more than $n$ vertices  such that $L(H)=L_{(m, n)}^T$.
\end{definition}

We will need the following lemma in the proof of our main results.

\begin{lemma} \label{lemma1} \cite{PZ} 
 There exists a left-compressed  extremal $T$-hypergraph $H$ for $L_{(m, n)}^T$.
\end{lemma}
{\bf Proof.} Let $H=(V, E)$ be an extremal $r$-graph for $\lambda_{(m, n)}^T$. Let ${\vec x}=(x_1, x_2, \ldots, x_n)$ be an optimal weight of $H$. We can assume that $x_i\ge x_j$ when $i<j$ since otherwise we can just relabel the vertices of $H'$ and obtain another  extremal $r$-graph with an optimal weight ${\vec x}=(x_1, x_2, \ldots, x_n)$ satisfying $x_i\ge x_j$ when $i<j$. If $H$ is not left compressed, performing a sequence of left-compressing operations (i.e. replace $E$ by $\mathcal{C}_{ij}(E)$ if $\mathcal{C}_{ij}(E)\neq E$), we will get a left-compressed $r$-graph $H'$ with the same number of edges, the same number of vertices, and $\lambda(H')\ge\lambda(H)$ . So $H'$ is a left-compressed extremal $r$-graph for $L_{(m, n)}^T$. \qed

\section{Proof of Theorem \ref{th2}}\label{proof-th2}

Let $H$ be a $\{1,2\}$- graph on $[n]$. In this case,
\begin{eqnarray*}
L_{\{\alpha_1 \}} (H,\vec{x})=\sum_{i \in E(H^1)} x_i+\alpha_1 \sum_{ i_1i_2  \in E(H^{2})}x_{i_1}x_{i_2},
\end{eqnarray*}
and
\begin{eqnarray}\label{L12}
 L_{\{\alpha_1 \}} (H) = \max \{ L (H, \vec{x}): \vec{x} \in S \}.
\end{eqnarray}
In the proof, we simply write $L_{\{\alpha_1 \}} (H,\vec{x})$ and $ L_{\{\alpha_1 \}} (H)$ as $L (H,\vec{x})$ and $ L (H)$.


\noindent{\bf \em  Proof of Theorem \ref{th2}.} Applying Lemma \ref{Lemmkkt} (a) and a direct calculation, we get a solution $\vec{y}$ to $(\ref{L12})$ when $H={K_t}^{\{ 1,2\} }$, which is given by $y_i=1/t$ for each $i$, $1\le i\le t$. Then $L(K^{\{1,2\}}_t)=1+\frac{\alpha_1}{2} - \frac{\alpha_1}{2t}$. Since $K^{\{1,2\}}_t \subseteq H$, then $L(H) \ge L(K^{\{1,2\}}_t).$
\\Now we proceed to show that $ L(H)\le  L(K^{ \{ 1,2 \} }_t)=1+\frac{\alpha_1}{2} - \frac{\alpha_1}{2t}$. Let ${\vec x}=(x_1, x_2, \ldots, x_n)$ be a solution to (\ref{L12}) satisfying (*) with $k$ positive weights. Without loss of generality, we may assume that $x_1 \ge x_2 \ge \ldots \ge x_k >x_{k+1}=x_{k+2}=\ldots = x_n=0$. By Lemma \ref{Lemmkkt} (b), $\forall 1\le i<j\le k,ij\in E(H^{2})$.

\begin{claim}
 $1$ $\forall 1\le i <j\le k$, if $i\in E(H)$ but $j\notin E(H)$, then $x_i-x_j=\frac{1}{\alpha_1}$.
\end{claim}
\noindent{\em Proof of Claim $1$.} By Lemma \ref{Lemmkkt} (a), ${\partial L(H, {\vec x})\over \partial x_i}={\partial  L(H, {\vec x})\over \partial x_j}$. By Lemma \ref{Lemmkkt} (b), $\forall 1\le i<j\le k,ij\in E(H^{2})$, therefore $1+\alpha_1(1-x_i)=\alpha_1(1-x_j)$, i.e. $x_i-x_j=\frac{1}{\alpha_1}$.  \qed

Let's continue the proof of Theorem \ref{th2}. Let $p=\lfloor \alpha_1 \rfloor$. Assume that there are $q$ 1-sets of $\{1, 2, 3, \cdots, k\}$ in $H^1$. If $q=k$, then $i\in E(H^1)$ for all $1 \le i\le k$, then $K^{\{1,2\}}_k$ is a subgraph of $H$. Since $t$ is the order of the maximum complete $\{1,2\}$-graph of $H$, then $k\le t$. We have
$$ L(H, {\vec x})=  L(K^{\{1,2\}}_k) = 1+\frac{\alpha_1}{2} - \frac{\alpha_1}{2k} \le 1+\frac{\alpha_1}{2} - \frac{\alpha_1}{2t}.$$

Therefore we can assume that  $q\le k-1$. Without loss of generality, assume that $i \in E(H^1)$ for $1\le i \le q $ and $j \notin E(H^1)$ for $q+1\le j\le k$,  by Claim $1$, $x_i=x_j+\frac{1}{\alpha_1}$, $\forall 1\le i\le q$ and $q+1\le j\le k$.
\\Case 1. $0 < \alpha_1 < 1$.
\\Note that $q=0$. Otherwise, by Claim $1$, we have $x_1={1 \over \alpha_1}+x_k>1$, which is a contradiction. So $\forall 1\le i\le k$, $i \not \in  E(H^1)$. Therefore,
$$L(H, {\vec x}) = L(K^{\{2\}}_{k})={\alpha_1 \over 2 }-{\alpha_1 \over 2k }\le 1+\frac{\alpha_1}{2} - \frac{\alpha_1}{2t}.
$$
\\Case 2. $\alpha_1 \ge 1$.
\\Note that $q \le p$. Otherwise, $x_1=x_k+\frac{1}{\alpha_1}$ , \ldots, $x_{q}=x_k+\frac{1}{\alpha_1}$, then $x_1+\ldots+x_{q}=qx_k+ \frac{q}{\alpha_1} > 1$, contradicts to $\sum_{i=1}^{k} x_i =1$ and $ x_i
> 0 {\rm \ for \ } 1\le i\le k$. Then $x_1= \ldots =x_q={1-q/\alpha_1 \over k}+1/\alpha_1$, $x_{q+1}=\ldots=x_k={1-q/\alpha_1 \over k}$,
\begin{eqnarray*}
 L(H, {\vec x})&=& x_1+\ldots+x_q+\alpha_1 \sum_{1\le i<j\le k}x_ix_j\\
&=& qx_1+\alpha_1 \sum_{1 \le i<j\le q} x_ix_j+\alpha_1 \sum_{1 \le i\le q}\sum_{q+1\le j \le k} x_ix_j+ \alpha_1 \sum_{q+1\le i<j\le k}x_ix_j\\
&=& qx_1 + \alpha_1{q \choose 2} x^2_1 + \alpha_1 q(k-q)x_1x_k+\alpha_1 {k-q\choose 2} x^2_k\\
&=& {q+\alpha_1 ^2 \over 2\alpha_1}-{ (\alpha_1 - q)^2 \over 2\alpha_1 k}.
\end{eqnarray*}
Next, we show ${q+\alpha_1 ^2 \over 2\alpha_1}-{ (\alpha_1 -q)^2 \over 2\alpha_1 k} < 1+\frac{\alpha_1}{2} - \frac{\alpha_1}{2t}$.
\begin{eqnarray*}
[1+\frac{\alpha_1}{2} - \frac{\alpha_1}{2t}]-[{q+\alpha_1 ^2 \over 2\alpha_1}-{ (\alpha_1 - q)^2 \over 2\alpha_1 k}]&=& 1-{q \over 2\alpha_1}+{ (\alpha_1 - q)^2 \over 2\alpha_1 k}- \frac{\alpha_1}{2t}\\
&>& 1-{q \over 2\alpha_1}- \frac{\alpha_1}{2t}.
\end{eqnarray*}
Since $p < \alpha_1 \le p+1$, $q \le p$  and $t \ge \alpha_1$, then ${q \over 2\alpha_1}+ \frac{\alpha_1}{2t}<1$. This completes the proof.
\qed

\section{Proof of Theorem \ref{thm4}}\label{proof-thm4}
Let $H$ be  a $\{1,r\}$- graph  with vertex set $V(H)=[n]$ and
edge set $E(H)$. In this case,
\begin{eqnarray*}
L_{\{\alpha_1 \}} (H,\vec{x})&=& \sum_{i \in E(H^1)} x_i+\alpha_1 \sum_{ i_1i_2\cdots i_r  \in E(H^{r})}x_{i_1}x_{i_2}\cdots x_{i_r},
\end{eqnarray*}
and
\begin{eqnarray} \label{L1r}
L_{\{\alpha_1 \}}(H)=max\{L(H,\vec x):\vec x \in S\}.
\end{eqnarray}
In the proof, we simply write $L_{\{\alpha_1 \}} (H,\vec{x})$ and $ L_{\{\alpha_1 \}} (H)$ as $L (H,\vec{x})$ and $ L (H)$.

\noindent{\bf \em Proof of Theorem \ref{thm4}.}
 Applying Lemma \ref{Lemmkkt} (a) and a direct calculation, we get a solution $\vec{y}$ to $(\ref{L1r})$ when $H={K_t}^{\{ 1,r\} }$ which is given by $y_i=1/t$ for each $i(1\le i\le t)$ and $y_i=0$ else. So $ L({{K_t}^{\{ 1,r\} }}) ={1+ \alpha_1 \frac{\prod_{i=1}^{r-1} (t-i)}{r!t^{r-1}}}$.
Since ${K_t}^{\{ 1,r\}}\subseteq H$, then $ L (H) \geq  L \left( {{K_t}^{\{ 1,r\} }}\right)$. Now we need to prove that $ L (H) \leq  L \left( {{K_t}^{\{ 1,r\} }}\right)$. Denote $ M_{\{t,\{1,r\}\}} =max\{ L (H):$ $H$ is a $\{1,r\}$-graph, $H$ contains a maximum complete subgraph $K_t^{\{1,r\}}$ and a maximum complete subgraph $K_t^{\{1\}}\}$.  If $ M_{\{t,\{1,r\}\}}\leq   L \left( {{K_t}^{\{ 1,r\} }}
\right) $, then $ L(H)\leq   L \left( {{K_t}^{\{ 1,r\} }} \right)$. Hence we can assume that $H$ is an extremal hypergraph, i.e., $ L(H)= M_{\{t,\{1, r\}\}}$.  If $H$ is not left-compressed, performing a sequence of left-compressing operations (i.e. replace $E$ by $\mathcal{C}_{ij}(E)$ if $\mathcal{C}_{ij}(E)\neq E$), we will get a left-compressed $\{1, r\}$-graph $H'$ with the same number of edges. The condition that the order of a maximum complete $\{1\}$-subgraph of $H$ is $t$ guarantees that
both the order of a maximum  complete $\{1, r\}$-subgraph of $H'$  and the order of a maximum complete $\{1\}$-subgraph of $H'$  are still $t$. By Lemma \ref{lem1}, $H'$ is an extremal graph as well.  So we can assume that the edge set of $H$ is left-compressed, $H^1=[t]$ and $[t]^{(r)}\subseteq H^r$. Let $\vec x = (x_1,\cdots, x_n)$ be a solution to (\ref{L1r}), where $x_1\geq x_2 \geq\ldots \geq x_k > x_{k+1}=x_{k+2}=\ldots =x_n=0$. If $k\le t$, then $ L(H)\le  L([k]^{\{1, r\}})\le  L([t]^{\{1, r\}})$.  So it suffices to show that $x_{t+1}=0$.

Let $1\le i\le t$. If $x_{t+1}>0$, then by Lemma \ref{Lemmkkt}, there exists $e\in E(H^r)$ such that $\{i, t+1\}\subset e$ and $\frac{{\partial  L \left( {H,\vec x}\right)}}{{\partial {x_i}}}=\frac{{\partial  L \left( {H,\vec x}\right)}}{{\partial {x_{t+1}}}}$.

Recall that $i \in E(H^1)$ and $t+1\notin E(H^1)$, then,
\begin{eqnarray*}
0 &=& \frac{{\partial  L \left( {H,\vec x} \right)}}{{\partial {x_i}}} -\frac{{\partial  L \left( {H,\vec x} \right)}}{{\partial {x_{t+1}}}}\\
&=& 1 +  L \left( {E^r_{i\backslash (t+1)} ,\vec x} \right)+ x_{t+1} L \left( {E^r_{i(t+1)} ,\vec x} \right)-x_i L \left( {E^r_{i(t+1)} ,\vec x} \right).
\end{eqnarray*}
Let $A= L \left( {E^r_{i(t+1)} ,\vec x} \right)$.
Thus, $x_i\geq \frac{1}{A}+x_{t+1}$. Note that $E^r_{i(t+1)}$ is a $(r-2)$-graph on $[n]\backslash \{i,t+1\}$, so $0<A \le \alpha_1 {(1-x_i-x_{t+1})^{r-2} \over (r-2)!}$. Then
\begin{equation}\label{eq2}
x_i> \frac{(r-2)!}{\alpha_1(1-x_i-x_{t+1})^{r-2}}+x_{t+1}.
\end{equation}

The above inequality clearly  implies that $x_i>{(r-2)! \over \alpha_1}$. If $ \alpha_1 \le (r-2)!$, then $x_i>1$ which is a contradiction. So what left is to consider $\alpha_1 > (r-2)!$.
Combining $x_i>{(r-2)! \over \alpha_1}$ with (\ref{eq2}), we have
\begin{equation}\label{eq1}
x_i> {(r-2)!\alpha_1^{r-3}\over [\alpha_1-(r-2)!]^{r-2}}.
\end{equation}
Recall that $\displaystyle {t\geq \lceil {[\alpha_1-(r-2)!]^{r-2} \over (r-2)!\alpha_1^{r-3}}\rceil }$, with the aid of (\ref{eq1}),
$\sum\limits_{i = 1}^t {{x_i}}  > 1 $,
a contradiction.
So $x_{t+1}=0$.
The proof is thus complete.
\qed

\section{Proof of Theorem \ref{th4}}\label{proof-th4}

Let $H$ be  a $\{1,2,3\}$- graph  with vertex set $V(H)=[n]$ and
edge set $E(H)$.  In this case,
\begin{eqnarray*}
 L_{\{\alpha_1,\alpha_2\}} (H,\vec{x})&=& \sum_{i \in E(H^1)} x_i+\alpha_1 \sum_{ i_1i_2 \in E(H^{2})}x_{i_1}x_{i_2}+\alpha_{2} \sum_{ i_1 i_2i_3 \in E(H^3)}x_{i_1}x_{i_2}x_{i_3},
\end{eqnarray*}
and
\begin{eqnarray} \label{L123}
L_{\{\alpha_1,\alpha_2\}}(H)=max\{L(H,\vec x):\vec x \in S\}.
\end{eqnarray}
In the proof, we simply write $L_{\{\alpha_1,\alpha_2\}} (H,\vec{x})$ and $ L_{\{\alpha_1,\alpha_2\}} (H)$ as $L (H,\vec{x})$ and $ L (H)$.

\noindent {\bf \em Proof of  Theorem \ref{th4}. }The proof is similar to the proof of Theorem \ref{thm4}.  Applying Lemma \ref{Lemmkkt} (a) and a direct calculation, we get a solution $\vec{y}$ to $(\ref{L123})$ when $H={K_t}^{\{ 1,r\} }$ which is given by $y_i=1/t$ for each $i(1\le i\le t)$ and $y_i=0$ else. So $ L \left( {{K_t}^{\{ 1, 2, 3\} }} \right) ={1+ \alpha_1 \frac{t-1}{2t}+\alpha_2 \frac{(t-1)(t-2)}{6t^2}}$. Hence we only need to prove $ L (H) =  L \left({{K_t}^{\{ 1, 2, 3\} }} \right)$. Since ${K_t}^{\{ 1, 2, 3\}} \subseteq H$,
then, $ L (H) \geq  L \left( {{K_t}^{\{ 1, 2, 3\} }} \right)$. Thus, to prove Theorem \ref{th4}, it suffices to prove that $ L (H) \leq  L \left( {{K_t}^{\{ 1, 2, 3\} }}
\right)$. Denote $ M_{\{t,\{1, 2, 3\}\}} =max\{ L (H):$ $H$ is a $\{1, 2, 3\}$-graph, $H$ contains a maximum complete subgraph $K_t^{\{1, 2, 3\}}$ and a maximum $\{1\}$ complete subgraph $K_t^{\{1\}}\}$.  If $ M_{\{t,\{1, 2, 3\}\}}\leq   L \left( {{K_t}^{\{ 1, 2, 3\} }}\right) $, then $ L(H)\leq   L \left( {{K_t}^{\{ 1, 2, 3\} }} \right)$. Hence we can assume that $H$ is an extremal hypergraph, i.e., $ L(H)= M_{\{t,\{1, 2, 3\}\}}$.  If $H$ is not left-compressed, performing a sequence of left-compressing operations (i.e. replace $E$ by $\mathcal{C}_{ij}(E)$ if $\mathcal{C}_{ij}(E)\neq E$), we will get a left-compressed $\{1, 2, 3\}$-graph $H'$ with the same number of edges. The condition that the order of a maximum complete $\{1\}$-subgraph of $H$ is $t$ guarantees that both the order of a maximum  complete $\{1, 2, 3\}$-subgraph of $H'$  and the order of a maximum  complete $\{1\}$-subgraph of $H'$  are still $t$. By Lemma \ref{lem1}, $H'$ is an extremal graph as well.  So we can assume that the edge set of $H$ is left-compressed, $H^1=[t]$, $[t]^{(2)}\subseteq H^2$ and $[t]^{(3)}\subseteq H^3$. Let $\vec x =(x_1,\cdots, x_n)$ be a solution to (\ref{L123}), where $x_1\geq x_2 \geq\ldots \geq x_k > x_{k+1}=x_{k+2}=\ldots =x_n=0$. If $k\le t$, then $ L(H)\le  L([k]^{\{1, 2, 3\}})\le  L([t]^{\{1, 2, 3\}})$.  So it suffices to show that $x_{t+1}=0$.

Let $1\le i\le t$. If $x_{t+1}>0$, then by Lemma \ref{Lemmkkt}, there exists $e\in E(H)$ such that $\{i, t+1\}\subset e$ and $\frac{{\partial  L\left( {H,\vec x}\right)}}{{\partial {x_i}}}=\frac{{\partial  L\left( {H,\vec x}\right)}}{{\partial {x_{t+1}}}}$. Let $ L ( E^2_{i(t+1)} ,\vec x)=\alpha_1$, if $i(t+1)\in E(H^2)$, and let $ L ( E^2_{i(t+1)} ,\vec x)=0$, if $i(t+1)\notin E(H^2)$. Recall that $i \in E(H^1)$ and $t+1\notin E(H^1)$,  then,
\begin{eqnarray*}
0 &=& \frac{{\partial  L\left( {H,\vec x} \right)}}{{\partial {x_i}}}-\frac{{\partial  L\left( {H,\vec x} \right)}}{{\partial {x_{t+1}}}}\\
 &=& 1 + L \left( {E^2_{i\backslash (t+1)} ,\vec x} \right)+x_{t+1} L \left( {E^2_{i(t+1)} ,\vec x} \right)+ L \left( {E^3_{i\backslash (t+1)} ,\vec x} \right)\\
&+& x_{t+1} L \left( {E^3_{i(t+1)} ,\vec x} \right)-x_{i} L \left( {E^2_{i(t+1)} ,\vec x} \right)-x_i L \left( {E^3_{i(t+1)} ,\vec x} \right).
\end{eqnarray*}

Let
$A= L \left( {E^2_{i(t+1)} ,\vec x} \right)+ L \left( {E^3_{i(t+1)} ,\vec x} \right)$.
Then $x_i\geq
\frac{1}{A}+x_{t+1}$, with $0<A \le \alpha_1+\alpha_2(1-x_i-x_{t+1})$. Hence
\begin{equation}\label{eq22}
x_i> \frac{1}{\alpha_1+\alpha_2(1-x_i-x_{t+1})}+x_{t+1}.
\end{equation}

The above inequality clearly  implies that $x_i>{1 \over \alpha_1+\alpha_2}$. If $\alpha_1+\alpha_2 \le 1$, then $x_i>1$ which is a contradiction. Combining $x_i>{1 \over \alpha_1+\alpha_2}$ with (\ref{eq22}), we have
\begin{equation}\label{eq11}
x_i> {\alpha_1+\alpha_2 \over (\alpha_1+\alpha_2)^2-\alpha_2 }.
\end{equation}
Recall that $t\ge \lceil{ (\alpha_1+\alpha_2)^2-\alpha_2 \over \alpha_1+\alpha_2}\rceil $, with the aid of (\ref{eq11}),
$\sum\limits_{i = 1}^t {{x_i}}  > 1 $, a contradiction.
So $x_{t+1}=0$.
The proof is thus complete.  \qed

\section{Proof of Theorem \ref{1.2}}\label{proof-1.2}
Let $H$ be  a $\{1,r_1,r_2,\ldots,r_l\}$- graph, $1<r_1<r_2< \ldots <r_l$, with vertex set $V(H)=[n]$ and
edge set $E(H)$.  In this case,

\begin{eqnarray*}
 L_{\{\alpha_1,\ldots ,\alpha_l \}} (H,\vec{x})= \sum_{i \in E(H^1)} x_i+\alpha_1 \sum_{ i_1i_2 \ldots i_{r_1}  \in E(H^{r_1})}x_{i_1}x_{i_2}\ldots x_{i_{r_1}}+\ldots
\end{eqnarray*}
\begin{equation*}
+\alpha_{l} \sum_{ i_1 i_2\ldots i_{r_l} \in E(H^{r_l})}x_{i_1}x_{i_2}\ldots x_{i_{r_l}},
\end{equation*}

and
\begin{eqnarray}\label{L1rm}
L_{\{\alpha_1,\ldots ,\alpha_l \}}(H)=max\{L(H,\vec x):\vec x \in S\}.
\end{eqnarray}
In the proof, we simply write $L_{\{\alpha_1,\ldots ,\alpha_l \}} (H,\vec{x})$ and $ L_{\{\alpha_1,\ldots ,\alpha_l \}} (H)$ as $L (H,\vec{x})$ and $ L (H)$.

\noindent{\bf{{\em Proof of Theorem \ref{1.2}}.}} Let $\vec x=(x_1,x_2,\ldots,x_n)$ be a solution to (\ref{L1rm}) satisfying(*) with $k$ positive weights. we may assume that $x_1 \ge x_2 \ge \ldots \ge x_k >x_{k+1}=\ldots = x_n$. Clearly, $ L(H) \ge  L(H[V(H^1)]). $
\begin{claim} \label{c1}
$2$ $(a)$ For all $1\le i \le k$, $i \in E(H^1) $.\\
$(b)$ Ether $i$ is isolated for all $1\le i \le k$ in $H[[k]]$ or $i$ is not isolated for any $1\le i \le k$ in $H[[k]]$.
\end{claim}
\noindent{\em  Proof of Claim $2$.} First we prove $(a)$. If $(a)$ fails to hold, then there are two possibilities.
\\Case 1. If $i \not \in E(H^1)$, for all $1\le i \le k$, then
\begin{eqnarray*}
 L(H,\vec x) &=& \alpha_1 \sum_{ i_1i_2 \ldots i_{r_1}  \in E(H^{r_1})}x_{i_1}x_{i_2}\ldots x_{i_{r_1}}+\ldots+\alpha_{l} \sum_{ i_1 i_2\ldots i_{r_l} \in E(H^{r_l})}x_{i_1}x_{i_2}\ldots x_{i_{r_l}}\\
&\le& \sum_{i=1}^l{\alpha_i \over r_i!}<\sum_{i=1}^l{\alpha_i \over (r_i-1)!}\le 1.
\end{eqnarray*}

But $ L(H,\vec x)= L(H) \ge  L(H^1)=1> L(H,\vec x)$, Contradiction!
\\Case 2. If there exists $1 \le i \le k$ and $1 \le j \le k$ such that $i \in E(H^1)$ and $j \not \in E(H^1)$, then
\begin{eqnarray*}
{\partial L (H, {\vec x})\over \partial x_i}&=&1+\alpha_1 \sum_{ i_1i_2 \ldots i_{r_1-1}  \in E_{i}^{r_1}}x_{i_1}x_{i_2}\ldots x_{i_{r_1-1}}+\ldots \\
& &+\alpha_{l} \sum_{ i_1 i_2\ldots i_{r_l-1} \in E_{i}^{r_l}}x_{i_1}x_{i_2}\ldots x_{i_{r_l-1}}>1,
\end{eqnarray*}
and
\begin{eqnarray*}
{\partial L (H, {\vec x})\over \partial x_j} &=& \alpha_1 \sum_{ i_1i_2 \ldots i_{r_1-1}  \in E_{j}^{r_1}}x_{i_1}x_{i_2}\ldots x_{i_{r_1-1}}+\ldots\\
& &+\alpha_{l} \sum_{ i_1 i_2\ldots i_{r_l-1} \in E_{j}^{r_l}}x_{i_1}x_{i_2}\ldots x_{i_{r_l-1}} < \sum_{i=1,\ldots,l}{\alpha_i \over (r_i-1)!} \le 1.
\end{eqnarray*}
${\partial L (H, {\vec x})\over \partial x_i}>{\partial L (H, {\vec x})\over \partial x_j}$, contradiction to  Lemma \ref{Lemmkkt} (a). So $(a)$ holds.
\\By $(a)$, we have $ L(H) \le  L(H[[k]])\le  L(H[V(H^1)])$, so $ L(H)= L(H[[k]]) =  L(H[V(H^1)])$.

Next,we prove $(b)$. By $(a)$, $i \in E(H^1)$ for all $1\le i \le k$. If there exists $1 \le i \le k$ and $1 \le j \le k$ such that $i$ is not isolated in $H[[k]]$ and $j$ is isolated in $H[[k]]$. Then there is some edge in $E(H[[k]])\backslash \{1,2,\ldots,k\}$ containing $i$ in but no  edge in $E(H[[k]])\backslash \{1,2,\ldots,k\}$ containing $j$, so
\begin{eqnarray*}
{\partial L(H, {\vec x})\over \partial x_i}=1&+&\alpha_1 \sum_{i_1i_2 \ldots i_{r_1-1} \in E_i^{r_1}}x_{i_1}x_{i_2}\ldots x_{i_{r_1-1}} \\
&+&\alpha_2 \sum_{i_1i_2 \ldots i_{r_2-1}  \in E_i^{r_2}}x_{i_1}x_{i_2}\ldots x_{i_{r_2-1}} \\
&+& \ldots + \alpha_l \sum_{i_1i_2 \ldots i_{r_l-1} \in E_i^{r_l}}x_{i_1}x_{i_2}\ldots x_{i_{r_l-1}} >1,
\end{eqnarray*}
but
\begin{eqnarray*}
{\partial L(H, {\vec x})\over \partial x_j}=1.
\end{eqnarray*}
Contradiction to ${\partial L(H, {\vec x})\over \partial x_i}={\partial L(H, {\vec x})\over \partial x_j}$. \qed

Let's continue the proof of Theorem \ref{1.2}, if $i$ is isolated for all $1\le i \le k$ in $H[[k]]$, then $ L(H)= 1$. If $i$ is not isolated for any $1\le i \le k$ in $H[[k]]$, then $i$ is not isolated for any $1\le i \le k$ in $H[V(H^1)]$. So $ L(H)= L(H[[k]])= L(H[V(H^1)\setminus D(H[V(H^1)])])$.

\section{Proof of Theorem \ref{1connection}}

Let $1<r_1<\ldots<r_l$ be positive integers and  let $\alpha_i(i=1,\ldots,l)$ be positive constants satisfying $\sum_{i=1}^l {\alpha_i \over (r_i-1)!} \le 1$. Let $m$ and $t$ be positive integers satisfying $t+\sum_{i=1}^l {t \choose r_i}< m\le t+1+\sum_{i=1}^l {t+1 \choose r_i}$.  
Let $T=\{1,r_1,r_2,\ldots,r_l\}$ and $Q=\{r_1,r_2,\ldots,r_l\}$. Let $H=(V, E)$ be an extremal $T$-graph for $L_{(m, n)}^{T}$. By Lemma \ref{lemma1} (b), we can assume that $H$ is left-compressed.
Let  ${\vec x}=(x_1, x_2, \ldots, x_n)$ be an optimal weighting of $H$ and $k$ be the number of non-zero weights in ${\vec x}$. Then $x_1\ge x_2 \ge \cdots \ge x_k>x_{k+1}=\cdots=x_n=0$. Since ${\vec x}$ has only $k$ positive weights, we can assume that $G$ is on $[k]$.

It is sufficient to show that $L(H, {\vec x})\le L(C_{m,T})$.   If $k\le t$, then $L(H)\le L(K^T)\le L(C_{m,T})$ since the range of  $m$ guarantees that $K_t^T \subset C_{m,T}$. So we can assume that $k\ge t+1$.

We first show the following result.

\begin{lemma}\label{1complete}
$E^1=[k]^{(1)}$.
\end{lemma}
{\bf Proof of Lemma \ref{1complete}.} If the lemma does not hold, then there are two possible cases.
Case 1. For each $i$,  $1\le i \le k$, $i \not \in E(H^1)$.
In this case,
\begin{eqnarray*}
 L(H,\vec x) &=& \alpha_1 \sum_{ i_1i_2 \ldots i_{r_1}  \in E(H^{r_1})}x_{i_1}x_{i_2}\ldots x_{i_{r_1}}+\ldots+\alpha_{l} \sum_{ i_1 i_2\ldots i_{r_l} \in E(H^{r_l})}x_{i_1}x_{i_2}\ldots x_{i_{r_l}}.
\end{eqnarray*}
Note that  for each $j, 1\le j\le l$,
$$\sum_{ i_1 i_2\ldots i_{r_j} \in E(H^{r_j})}x_{i_1}x_{i_2}\ldots x_{i_{r_j}} \le  \sum_{ i_1 i_2\ldots i_{r_j} \in E([k]^{(r_j)})}x_{i_1}x_{i_2}\ldots x_{i_{r_j}}$$
 and $\sum_{ i_1 i_2\ldots i_{r_j} \in E([k]^{(r_j)})}x_{i_1}x_{i_2}\ldots x_{i_{r_j}}$ reaches the maximum ${{k \choose r_j} \over k^{r_j}}\le {1 \over r_j!}$ when $\vec x$ is the characteristic vector of $[k]$. Therefore,
$$L(H,\vec x)\le \sum_{i=1}^l{\alpha_i \over r_i!}<\sum_{i=1}^l{\alpha_i \over (r_i-1)!}\le 1.$$
But $ L(H,\vec x)= L(H) \ge  L(H^1)=1> L(H,\vec x)$. Contradiction!
\\Case 2. There exists $1 \le i \le k$ and $1 \le j \le k$ such that $i \in E(H^1)$ and $j \not \in E(H^1)$.
In this case,
\begin{eqnarray*}
{\partial L (H, {\vec x})\over \partial x_i}&=&1+\alpha_1 \sum_{ i_1i_2 \ldots i_{r_1-1}  \in E_{i}^{r_1}}x_{i_1}x_{i_2}\ldots x_{i_{r_1-1}}+\ldots \\
& &+\alpha_{l} \sum_{ i_1 i_2\ldots i_{r_l-1} \in E_{i}^{r_l}}x_{i_1}x_{i_2}\ldots x_{i_{r_l-1}}>1,
\end{eqnarray*}
and
\begin{eqnarray*}
{\partial L (H, {\vec x})\over \partial x_j} &=& \alpha_1 \sum_{ i_1i_2 \ldots i_{r_1-1}  \in E_{j}^{r_1}}x_{i_1}x_{i_2}\ldots x_{i_{r_1-1}}+\ldots\\
& &+\alpha_{l} \sum_{ i_1 i_2\ldots i_{r_l-1} \in E_{j}^{r_l}}x_{i_1}x_{i_2}\ldots x_{i_{r_l-1}} < \sum_{i=1,\ldots,l}{\alpha_i \over (r_i-1)!} \le 1.
\end{eqnarray*}
So ${\partial L (H, {\vec x})\over \partial x_i}>{\partial L (H, {\vec x})\over \partial x_j}$, contradiction to  Lemma \ref{Lemmkkt}. 
This completes the proof of Lemma \ref{1complete}. \qed

Let us continue the proof of the theorem. By Lemma \ref{1complete},  $L (H, {\vec x})=1+L (H^Q, {\vec x})$, and the number of edges in $H^Q$ is $m-k\le m-t-1$. By the assumption, $L (H^Q)\le L(C_{m-t-1, Q})$. Therefore, $L(H)=L (H, {\vec x}) \le 1+L(C_{m-t-1, Q})$. Note that $E(C_{m, T})=[t+1]^{(1)}\cup E(C_{m-t-1, Q})$. So $1+L(C_{m-t-1, Q})=L(C_{m, T})$. Consequently, $L(H)\le L(C_{m, T})$. This completes the proof of this theorem. \qed


\section{Conclusion}
The classical method of Lagrange multiplier has been applied often in evaluating the Graph-Lagrangian of a hypergraph. However, evaluating the Graph-Lagrangian of a general hypergraph seems to be  challenging  and very few general results are known for hypergraphs. In the future, we will learn and explore whether modern Lagrange theory \cite{Giannessi} will help advance the research.

\section{ Acknowledgments} 

We thank Franco Giannessi for helpful discussions.



\begin{thebibliography}{1}

\bibitem{PPTZ}
Y. Peng, H. Peng, Q. Tang, C. Zhao, An extension of Motzkin-Straus
Thorem to non-uniform hypergraphs and its applications, preprint, http://arxiv.org/abs/1312.4135v1.

\bibitem{GLPS} Gu, R. , Li, X. , Peng, Y., Shi, Y.: Some Motzkin-Straus type results for non-uniform hypergraphs, preprint, http://arxiv.org/abs/1310.8442.

\bibitem{MS}
T. Motzkin, E. Straus, Maxima for graphs and a new proof of a
theorem of Tur\'{a}n, {\it Canad. J. Math.} {\bf 17}(1965),
533--540.

\bibitem{B1}  Bomze, I.M.: Evolution towards the maximum clique. J. Glob. Optim. \textbf{10},  143-164  (1997)

\bibitem{B2} Budinich, M.: Exact bounds on the order of the maximum clique of a graph. Discret Appl. Math. \textbf{127},  535-543  (2003)

\bibitem{B3} Busygin, S.: A new trust region technique for the maximum weight clique problem. Discret Appl. Math. \textbf{154},  2080-2096  (2006)

\bibitem{G9} Gibbons, L.E., Hearn, D. W.,  Pardalos, P. M.,  Ramana, M. V.: Continuous characterizations of the maximum clique problem. Math. Oper. Res. \textbf{22},  754-768  (1997)

\bibitem{PP}P. M. Pardalos and A. T. Phillips, A global optimization approach for solving the maximum clique problem,  Int. J. Compt. Math. 33 (1990), 209-216.

\bibitem{Turan}
P. Tur\'an, On an extremal problem in graph theory (in Hungarian), {\it Mat. Fiz. Lapok}
{\bf 48} (1941), 436--452.

\bibitem{Keevash}
P. Keevash, Hypergrah Tur\'an problems, {\it Surveys in Combinatorics}, Cambridge
University Press, (2011), 83--140.

\bibitem{Sidorenko}
A. Sidorenko, The maximal number of edges in a homogeneous
hypergraph containing no prohibited subgraphs, {\it Math Notes} {\bf
41}(1987), 247--259. Translated from Mat. Zametki.

\bibitem{FF}
P. Frankl, Z. F\"uredi, Extremal problems and the Lagrange
function of hypergraphs, {\it Bulletin Institute Math. Academia
Sinica} {\bf 16}(1988), 305--313.

\bibitem{Mu}
D. Mubayi, A hypergraph extension of Turan¡¯s theorem, {\it J.
Combin. Theory Ser. B} {\bf 96}(2006), 122--134.


\bibitem{keevash2}
D. Hefetz and P. Keevash, A hypergraph Tur\'an theorem via lagrangians of intersecting
families, {\it J. Combin. Theory Ser. A} {\bf 120} (2013), 2020--2038.

\bibitem{talbot} J. Talbot, Lagrangians of hypergraphs, {\it Comb. Probab. Comput.} {\bf 11} (2002), 199--216.

\bibitem{TPZZ}Q. S. Tang, Y. Peng, X. D. Zhang, and C. Zhao, Lagrangians of $3$-uniform Hypergraphs, preprint.

\bibitem{TPZZ1} Q. S. Tang, Y. Peng, X. D. Zhang, and C. Zhao, Some results on Lagrangians of Hypergraphs, {\it Discrete Appl. Math} 166(2014), 222-238.

\bibitem{STZP}Y. Sun, Q. S. Tang, C. Zhao and Y. Peng, On the largest Graph-Lagrangian of 3-graphs with fixed number of edges, {\it J. Optimiz. Theory App.} (accepted).

\bibitem{GK}
J. Griggs, G. Katona, No four subsets forming an N, {\it J. Combin.
Theory Ser. A.} {\bf 115}(2008), 677--685.

\bibitem{GL}
J. Griggs, L. Lu, On families of subsets with a forbidden subposet,
{\it Comb. Probab. Comput.} {\bf 18}(2009), 731--748.

\bibitem{JL}
T. Johston, L. Lu, Tur\'an problems on non-uniform hypergraphs,
submitted.

\bibitem{FF1}
P. Frankl, Z. F\"uredi, Extremal problems whose solutions are the
blow-ups of the small Witt-designs, {\it J. Combin. Theory Ser. A.}
{\bf 52}(1989), 129--147.

\bibitem{FR}
P. Frankl, V. R\"{o}dl, Hypergraphs do not jump, {\it Combinatorica}
{\bf 4}(1984), 149--159.

\bibitem{PZ}
Y. Peng, C. Zhao, A Motzkin-Straus type result for $3$-uniform
hypergraphs, {\it Graphs Combin.} {\bf 29}(2013), 681--694.

\bibitem{Giannessi} 
F. Giannessi, "Constrained Optimization and Image Space Analysis. Vol.I: Separation of Sets and Optimality Conditions", Springer, New York, 2005, pp.1-395.

\bibitem{BP} S. R. Bul\'o and M. Pelillo, A generalization of the Motzkin-Straus theorem to hypergraph, Optim Lett 3(2009), 287-295.

\end{thebibliography}
\end{document}